\documentclass[12pt,letter]{article}
\usepackage{hyperref}
\usepackage{graphicx}
\usepackage{mathptmx}       
\usepackage{helvet}         
\usepackage{courier}        
\usepackage{amssymb}
\usepackage{tikz,amsmath,listings,textcomp}
\usepackage{xfrac,subfigure}
\usepackage{makeidx}         
\usepackage{graphicx}        
\usepackage{multicol}        
\usepackage[bottom]{footmisc}

\newcommand{\e}{\ensuremath{\varepsilon}}

\renewcommand{\emph}[1]{\textsl{#1}}
\usepackage{bm}

\title{Inverse Cubic Iteration}
\author{Robert M.~Corless}

\begin{document}
\maketitle
\begin{abstract}
  There are thousands of papers on rootfinding for nonlinear scalar equations. Here is one more, to
  talk about an apparently new method, which I call ``Inverse Cubic Iteration'' (ICI) in analogy to the Inverse Quadratic Iteration in Richard Brent's \texttt{zeroin} method. The possibly new method is based on a cubic blend of tangent-line approximations for the inverse function.  We rewrite this iteration for numerical stability as an \emph{average} of two Newton steps and a secant step: only one new function evaluation and derivative evaluation is needed for each step.  The total cost of the method is therefore only trivially more than Newton's method, and we will see that it has order $1+\sqrt{3} = 2.732...$, thus ensuring that to achieve a given accuracy it usually takes fewer steps than Newton's method while using essentially the same effort per step.
\end{abstract}

\section{Introduction}
Recently I wrote, with Erik Postma, a paper on \emph{blends} of Taylor series~\cite{Corless:2020:MapleBlends}.  A ``blend'' is just a two-point Hermite interpolant: if we have series for $f(x)$ at $x=a$ and another series for the same function $f(x)$ at $x=b$ then we can blend the series to get a (usually) better approximation to $f(x)$ across the interval.  The paper discusses how to compute high-order blends stably and efficiently in Maple, and looks at some applications of them.

In contrast, this paper looks at a low-order blend, and a particular application, namely rootfinding.  In fact all we will need is a cubic.  If we just use series up to and including terms of order $1$ at each point, namely $f(a)$ and $f'(a)$ and $f(b)$ and $f'(b)$ then these four pieces of information determine, as is well known, a cubic Hermite interpolant.

What is less well-known (but still perfectly well-known, to be honest) is that from these same pieces of information one can usually construct a cubic interpolant of the \emph{inverse} function inv$f(y)$, also written $\breve{f}(y)$. The data for that information is that the points of expansion are $y=f(a)$ and $y=f(b)$, and the function values are $a$ and $b$, and the derivative values are $1/f'(a)$ and $1/f'(b)$: four pieces of information which together suffice to determine a cubic polynomial in $y$ that fits the data.

Of course we cannot have $f'(a) = 0$ or $f'(b) = 0$, and indeed $f'(x)=0$ anywhere near the interval is bad news for a rootfinder: the root will be ill-conditioned, if it's determined at all.  If $f'(x)=0$ anywhere in the interval then the function is not guaranteed to be monotonic, and so the inverse function will not be defined in the whole neighbourhood.  Throwing caution to the winds, we're not going to worry about this at all.
\section{Derivation of the Formula}
Nearly any numerical analysis textbook will have an expression equivalent to the following for a cubic Hermite interpolant, fitting four pieces of information to get a cubic polynomial that fits those pieces of information:
\begin{align}\label{eq:cubichermite}
  y =& \left( 1+2\,\theta \right)  \left(\theta-1 \right) ^{2}f(a) +
\theta\, \left( \theta-1 \right) ^{2}{hf'(a)}\nonumber\\
&+{\theta}^{2}
 \left( 3-2\,\theta \right) f(b)+{\theta}^{2} \left( \theta-1
 \right) {h f'(b)}
\end{align}
where $\theta = (x-a)/(b-a)$ and $h= b-a$. When $\theta=0$, $x=a$; when $\theta=1$, $x=b$.  Since the derivative $d/d\theta = dx/d\theta \cdot d/dx = h d/dx$ we see that it all works out: $y = f(a)$ when $\theta=0$, $dy/d\theta = hf'(a)$ and $y=f(b)$ when $\theta=1$ and $dy/d\theta = hf'(b)$ when $\theta = 1$.

Now we do the inverse.  Put $s = (y - f(a))/(f(b)-f(a))$ and $\Delta = f(b)-f(a)$.  We see right away that we are in trouble if $f(a)=f(b)$; but then, by Rolle's theorem there would be a place $c$ between $a$ and $b$ where $f'(c)=0$ if that were the case, and we have already said that we are going to ignore that difficulty.  Indeed, both Inverse Quadratic Iteration (IQI) and the secant method \emph{also} fail in that case, so at least we are in good company.

We want
\begin{align}\label{eq:invcubichermite}
  x =& \left(  \left( 1+2\,s \right) a+{\frac {s\Delta}{f'\left( a \right) }} \right)  \left( 1-s \right) ^{2}\nonumber\\
 &+
 \left(  \left( 3-2\,s \right) b-{\frac { \left( 1-s \right) \Delta}{
 f' \left( b \right) }} \right) {s}^{2}
\>.
\end{align}
By similar reasoning to before, we see that at $s=0$ then $x=a$ and $dx/dy = 1/f'(a)$, while at $s=1$ then $x=b$ and $dx/dy = 1/f'(b)$.

Why do this?  The wonderfully simple trick of Inverse Quadratic Iteration, namely that to find an approximate value of $x$ that makes $y=0$ one simply substitutes $y=0$ into the Inverse formula, also works here. Putting $y=0$ means $s = (0-f(a))/(f(b)-f(a))$ (a number between $0$ and $1$) and that gets us our definite value of $x$. Doing this blindly, we get
\begin{align}
&\left(  \left( 1-2\,{\frac {f \left( a \right) }{\Delta}} \right) a-{
\frac {f \left( a \right) }{f'  \left( a
 \right) }} \right)  \left( 1+{\frac {f \left( a \right) }{\Delta}}
 \right) ^{2}\nonumber\\
 &\qquad+{\frac { \left( f \left( a \right)  \right) ^{2}}{{
\Delta}^{2}} \left(  \left( 3+2\,{\frac {f \left( a \right) }{\Delta}}
 \right) b-{\frac {\Delta}{f'  \left( b
 \right) } \left( 1+{\frac {f \left( a \right) }{\Delta}} \right) }
 \right) }
\end{align}
This is correct, but I suspect numerically useless.  To be honest, I didn't even try it.  Instead, I looked for a way to rewrite it to improve its numerical stability.

\section{A numerically more stable expression}
After a surprisingly small amount of algebraic manipulation (perhaps a generous share of luck was involved) I arrived at the following:
\begin{align}
   x_{n+1} =& {\frac {{y_{{n}}}^{2}}{ \left( y_{{n-1}}-y_{{n}} \right) ^{2}} \left(
x_{{n-1}}-{\frac {y_{{n-1}}}{f'  \left( x_{{n-
1}} \right) }} \right) }+{\frac {{y_{{n-1}}}^{2}}{ \left( y_{{n-1}}-y_
{{n}} \right) ^{2}} \left( x_{{n}}-{\frac {y_{{n}}}{f'  \left( x_{{n}} \right) }} \right) }\nonumber \\
            & -2\,{\frac {y_{{n-1}}y_{{
n}}}{ \left( y_{{n-1}}-y_{{n}} \right) ^{2}} \left( x_{{n}}-{\frac {y_
{{n}} \left( x_{{n}}-x_{{n-1}} \right) }{y_{{n}}-y_{{n-1}}}} \right) }\>. \label{eq:ICIiteration}
\end{align}
Call this iteration Inverse Cubic Iteration, or ICI.
Instead of presenting a tedious transformation of the previous formula into the one above (it is equivalent, trust me, although I did make the sleight-of-hand transformation from $a$ and $b$ to $x_{n-1}$ and $x_n$, implying iteration), I will verify that it has the appropriate behaviour in the next section; we will look there at how quickly it converges to a root (when it converges at all).
The interesting part of this, to me, is that we can \emph{recognize} all three terms in equation~\eqref{eq:ICIiteration}.  There is a Newton iteration starting at $x_{n-1}$, another at $x_n$, and a \emph{secant} iteration using both points.  Then all three are averaged together\footnote{It might be numerically better still to average the previous points $a_0x_n + a_1x_{n-1} + a_sx_{n}$, and average the small updates $a_0y_n/f'(x_n) + a_1 y_{n-1}/f'(x_{n-1}) + a_s(y_n\Delta x/\Delta y)$, and then use this average update to improve the average of the previous points.  Currently this is the way my implementation does it, but it may not make much difference.}: the weights are $y_n^2$, $y_{n-1}^2$, and $-2y_ny_{n-1}$ all divided by $(y_n-y_{n-1})^2$ so that the weights all add up to $1$, as they should.  The most accurate estimate is given the most weight, by using the residuals $y_k$ on the other estimates.  I think this is a rather pretty formula.

The numerical stability (or lack thereof) of the secant method and of Newton's method is well-understood; by writing the secant method as above we make a ``small change'' to an existing value, and similarly Newton's method is written in as stable a way as it can be. Thus we have what is likely to be a reasonably stable iteration formula.

Because this method uses two points $x_{n-1}$ and $x_n$ to generate a third point $x_{n+1}$, it seems similar to the secant method.  Because it uses derivatives $f'(x_k)$ at each point, it seems similar to Newton's method.  It is not very similar to IQI, which fits an inverse quadratic to three points and uses that to generate a fourth.  Halley's method and other Schr\"oder iterations use higher derivatives and will be more expensive per step.  [It might be interesting to try to average higher-order methods in a similar manner to speed up convergence even more; so far I have not been able to succeed in doing so.  The simple trick of replacing $f$ by $f/\sqrt{f'}$ which transforms Newton's method to Halley's method (and speeds up secant) does not seem to work, with these weights.]

I am going to assume that we start with a single initial guess $x_0$ and then generate $x_1$ by a single step of Newton's method.  Then it is obvious that to carry out one step with equation~\eqref{eq:ICIiteration} we need only \emph{one} more function evaluation $y_1 = f(x_1)$ and \emph{one} more derivative evaluation $f'(x_1)$, because we can record the previous $y_0 = f(x_0)$ and $f'(x_0)$ instead of recomputing them.  In the standard model of rootfinding iteration where the dominant cost is function evaluation and derivative evaluation, we see that this method is no more expensive than Newton's method itself.  But instead of computing $X_{n+1} = x_n - f(x_n)/f'(x_n)$ and forgetting all about $x_{n-1}$, once we have formed this Newton iteration we average it with the previous Newton iteration together with a secant iteration: the cost of forming this average is no more than a few floating-point operations, which we may consider to be trivial compared to the assumed treasury-breaking expense of evaluating $f(x)$ and $f'(x)$. When we try this on Newton's classical example $z^3-2z-5$ with an initial guess of $z_0 = 1$ we get $29$ Digits of accuracy after $7$ iterations ($10$ Digits of accuracy after $6$ iterations).

I will sketch the cost of this method using the normal conventions in section~\ref{sec:cost}.
\section{The order of the method}
If we assume that $r$ is the root we are searching for, $f(r) = 0$, and that $x_{n-1} = r + \e_{n-1}$ and $x_n = r + \e_{n}$ then a straightforward and brutal series computation in Maple shows that
\begin{equation}\label{eq:errorformula}
  x_{n+1} = r + (\e_{n-1} \e_{n})^2\frac{f^{(iv)}(r)f''(r) - 10 f'(r)f''(r)f'''(r) + 15 (f''(r))^3 }{4!(f'(r))^3} + \cdots\>.
\end{equation}
The error in the residual, which is what we will actually measure, is ever so slightly different: $y(x_k) = f(r + \e_k) = f(r) + f'(r)\e_k$ and so we can remove one $f'(r)$ factor from the denominator term above.  Notice that each $x_k$ exactly solves $F(x) = f(x) - y_k$; this is a trivial observation, but often useful: by iterating, we get the exact solutions to slightly different equations.  Here, though, it doesn't matter.  Both formulae show the trend of the errors.
That is, the error in the next iterate is the \emph{square of the product} of the two previous errors.  This is satisfyingly analogous to that of the secant method (product of the two previous errors) and Newton's method (square of the previous error).  Solving the recurrence relation $\ln\e_{n+1} = 2(\ln\e_n + \ln\e_{n-1})$ shows that the error is asymptotically the $1+\sqrt{3}=2.732\ldots$th power of the previous one; not quite cubic convergence, but substantially faster than quadratic convergence.  

\begin{figure}
  \centering
  \includegraphics[width=8cm]{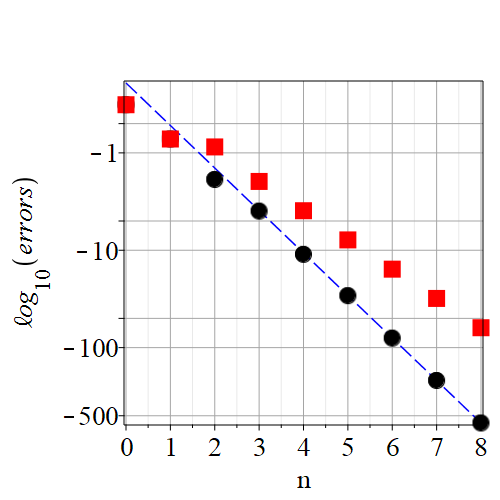}
  \caption{Convergence to a root of $f(x) = (x^2+x)\exp(-x) - 1/3$.  The initial guess was $x_0 = 2.0$. The next iterate was computed by Newton's method, $x_1 = x_0 - f(x_0)/f'(x_0)$.  Thereafter equation~\eqref{eq:ICIiteration} was used.  Computing was carried out in $1000$ Digits.  What is plotted is the logarithm (base $10$) of the residual error, equivalently $\log_{10}|y_n| = \log_{10}|f(x_n)|$. The dashed blue line represents the theoretical curve $y_n = (0.6437)^{(1+\sqrt{3})^n}$ demonstrating the theoretical computation of the order of the method above. The constant $0.6437$ was found by fitting the final data point. Newton's method (red squares) achieves only $63$ Digits of accuracy instead of ICI's nearly $600$. }\label{fig:ICIexp}
\end{figure}

Let's take the method out for a spin to see if practice matches theory.  Consider the function $f(x) = (x^2+1)\exp(-x) - 1/3$.  This has two places where $f'(x)=0$, namely $x = (1 \pm \sqrt{5})/2$.  Provided we stay away from there, we should be fine.  Taking our initial guess to be $x_0 = 2.0$, we then compute $x_1$ by Newton's method as usual.  Now we have two iterates to work with; we then compute $x_2$, $x_3$, $\ldots$, $x_8$.  Because convergence is so rapid, and I wanted to demonstrate it, I worked with Digits set to be $1000$.  This was overkill, but it allows the plot in figure~\ref{fig:ICIexp} to show clearly that the error in $x_8$ is about $10^{-594}$.  Newton's method (red squares) achieves only $10^{-63}$ with the same effort. The ratios $y_k/(y_{k-1}y_{k-2})^2$ for $k=2$, $3$, $\ldots$, $8$ were, respectively, $1.5952$, $17.048$, $4.5955$, $4.9061$, $4.9080$, $4.9081$, and $4.9080$. These fit the curve $y_k = (0.6437)^{(1+\sqrt{3})^k}$, demonstrating the near-cubic order of the method. The constant $0.6437$ was found by fitting to the final data point.  This trend predicts that the accuracy in $x_9$ would be $ 4.9081(y_8y_7)^2 = 1.7383\cdot 10^{-1622}$.  Re-doing the computation in $1624$ Digits and extending to $x_9$, we get $y_9 = 1.7383\cdot 10^{-1622}$, just as predicted.
\section{The cost of the method\label{sec:cost}}
The best analysis---that I know---of the cost of Newton's method for solving a scalar equation by iterating with constant precision is contained in~\cite{Neumaier(2001)}.  There, the author weighs the cost of function evaluations and estimates that the cost of derivative evaluations are usually a modest factor larger than that of function evaluations; then, balancing the greater cost of Newton's method per step against its general need for fewer iterations, he concludes that usually the secant method should win.  Yes, the secant method usually takes more iterations; but each one usually costs only a fraction of a Newton step.  It is \emph{only} when the derivatives are unusually inexpensive that Newton's method wins.  For the example in figure~\ref{fig:ICIexp} this is the case: the dominant cost of the function evaluation is the exponential, and that can be re-used in the derivative.  So, for this particular problem, ICI should win over Newton's method and even more decisively over the secant method. But to win over Newton's method, ICI needs to take fewer iterations, because the cost per step is essentially identical to Newton's method.  ICI will never\footnote{``What, never?''  ``No, never!'' ``What, \emph{never}?'' ``Hardly ever!'' ---the crew and Captain of HMS Pinafore} take \emph{more} iterations than Newton's method.  But will it take fewer?

That is actually the issue.  Will we save any iterations at all?  Newton's method typically converges very rapidly, given a good initial guess.  To give a very rough sketch, if the initial error is $2^{-2}$ and all the constants are $1$, we can expect a double precision answer in five iterations (actually the error then would be wasted: $2^{-64}$, and rounding errors would have made it at best $2^{-52}$).  ICI would have the sequence of errors $\e_0=2^{-2}$, $\e_1=2^{-4}$ (same as Newton for the $x_1$ iterate, of course), $\e_2=2^{-12}$ squaring the product $2^{-2}\cdot 2^{-4}$, then $\e_3 = 2^{-2(12+4)} = 2^{-32}$ which is the same error that Newton achieved in \emph{four} iterations.  One more gets us $\e_4 = 2^{-2(32+12)} = 2^{-88}$ which is a vast overkill; still, we have beaten Newton by one iteration with this initial error.

If instead the initial guess is only accurate to $2^{-1}$, not $2^{-2}$, then Newton takes six iterations while ICI \emph{also} takes six: $\e_0 = 2^{-1}$, $\e_1=2^{-2}$, $\e_3 = 2^{-6}$, $\e_4 = 2^{-16}$, $\e_5 = 2^{-44}$.

If instead the initial guess is accurate to $2^{-3}$ then Newton takes five iterations, while ICI takes $\e_1 = 2^{-6}$, $\e_2 = 2^{-18}$, $\e_3 = 2^{-48}$ and at four iterations beats Newton.

The lesson seems clear: if one wants ultra-high precision, the ICI method will give the same results in fewer iterations \emph{of essentially the same cost}.  If one only wants a double precision result, then the number of iterations you save may not be much, if any (and if it doesn't really beat Newton, then it won't often beat the secant method).  And \emph{both} methods demand good initial guesses.

\section{Initial guesses}
I have been using the scheme ``choose only one initial guess and let Newton's method determine the next iterate'' above.  Obviously there are alternative strategies.  I like this one because it's usually hard enough to come up with \emph{one} initial guess to a root, let alone two.

\begin{figure}
  \centering
  \includegraphics[width=6cm]{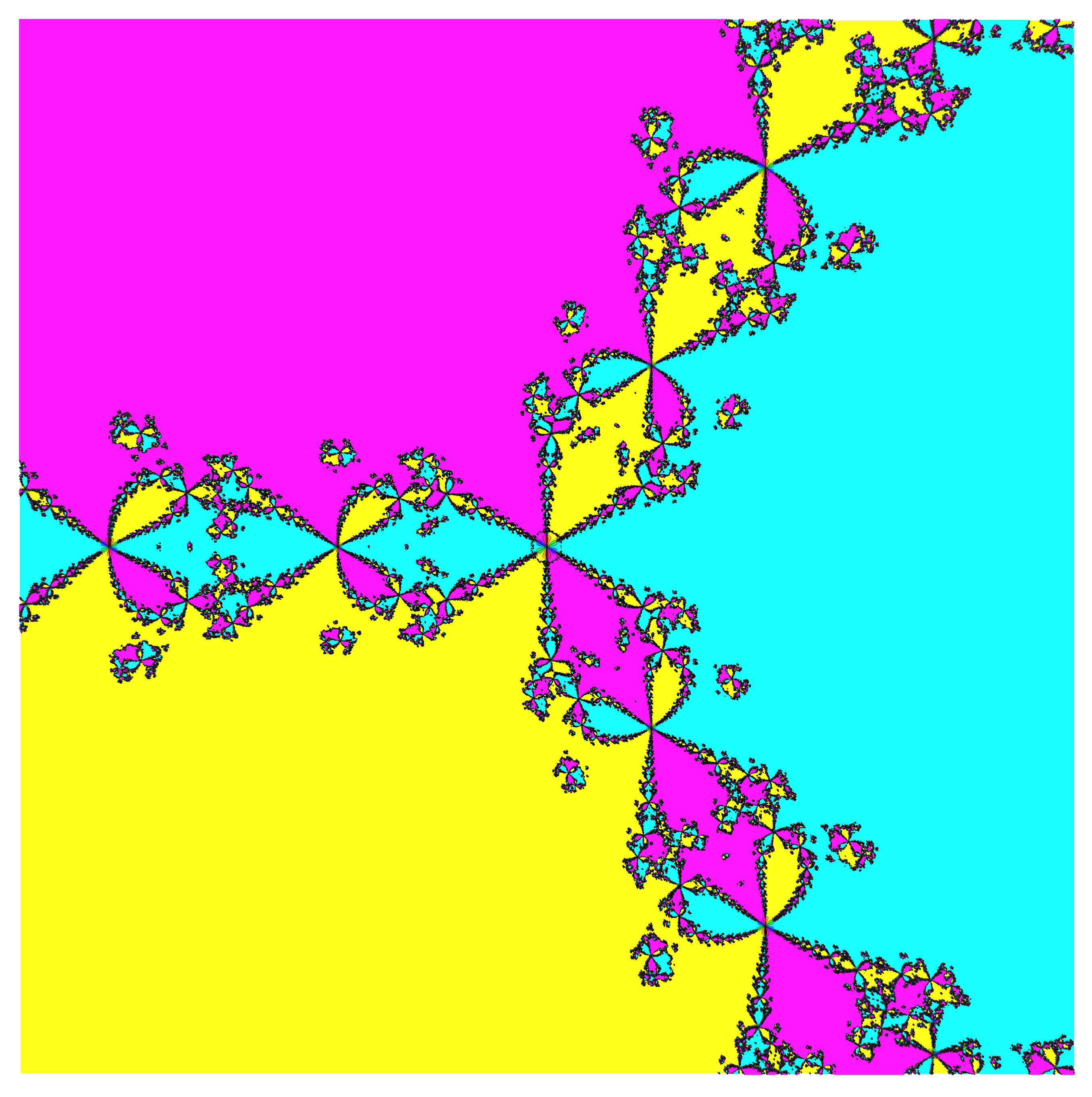}
  \caption{Approximate basins of attraction for the ICI method for $f(z) = z^3-1$. Each initial guess $z_0 = x+iy$ is coloured with the phase of the result after $13$ iterations of ICI, namely arg$(z_{13})$ (earlier iterations may have converged, in which case the limit is recorded).  The first iterate $z_1$ is obtained by Newton's method. Computation done on a $1600$ by $1600$ grid on $-2 \le x \le 2$, $-2 \le y \le 2$. Notice that the fractal structure is more complicated than that for simple Newton's method, suggesting that the choice of initial guess is more fraught for this method. }\label{fig:fractal}
\end{figure}

\begin{figure}
  \centering
  \includegraphics[width=6cm]{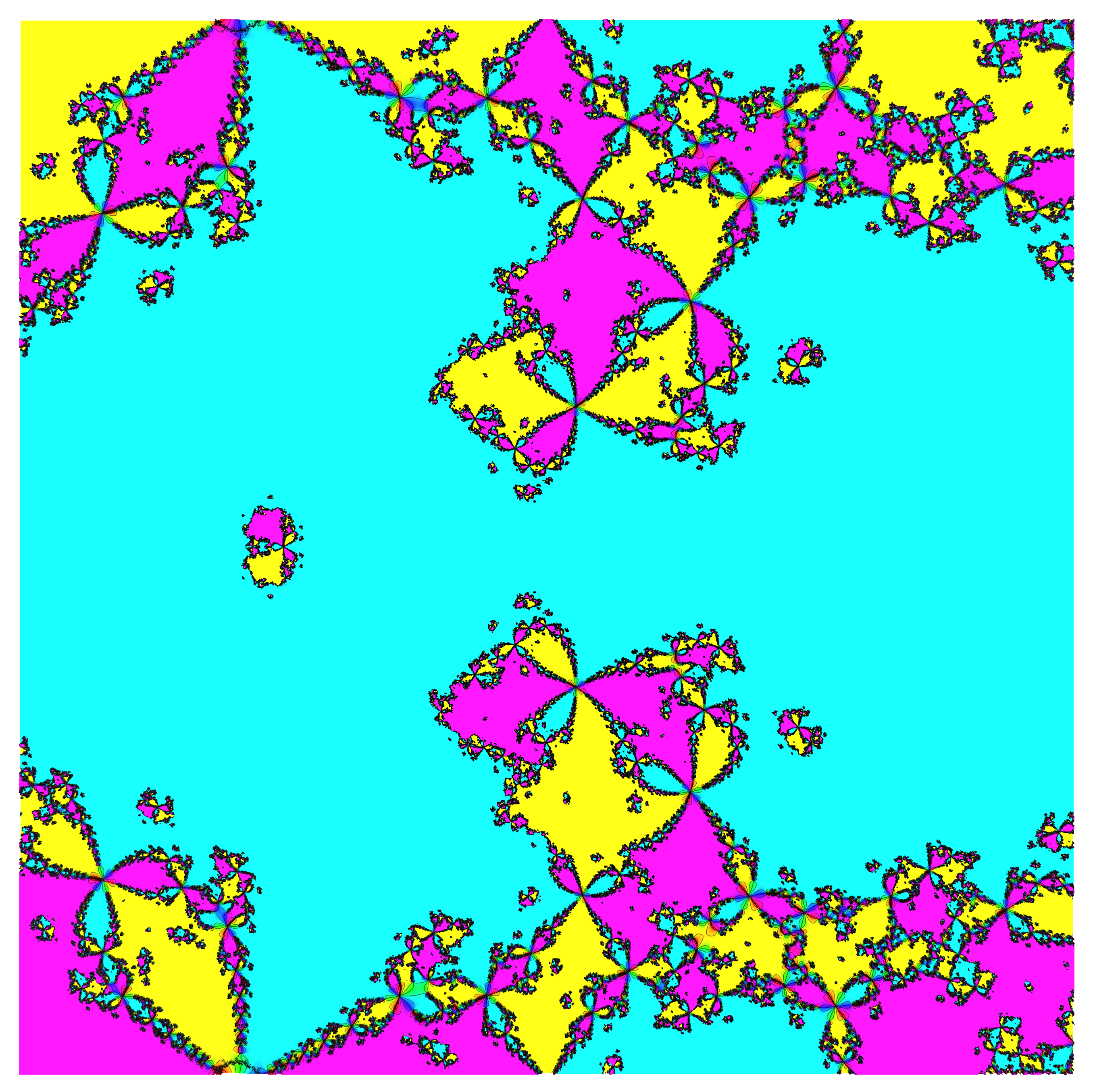}
  \caption{Zooming in, we show the region $-1.45 \le x \le -1.05$ and $-0.2 \le y \le 0.2$ again with a $1600$ by $1600$ grid for the same function and iterations as in figure~\ref{fig:fractal}.  We see disconnected components. }\label{fig:fractalzoom}
\end{figure}

Using this strategy also allows one to plot basins of attraction: which initial guesses converge to which root?  I took the function $f(z) = z^3-1$ and sampled a $1600$ by $1600$ grid on $-2 \le x \le 2$, $-2 \le y \le 2$ where $z= x+iy$ and took at most $13$ iterations of ICI starting from each gridpoint; if the iteration converged (to $10^{-8}$) earlier, that earlier result was recorded.  Plotting the argument (phase) of the answer gives an apparent fractal diagram, very similar to that for Newton's method, but different in detail.  In particular, the basins of attraction are visibly disconnected.  See figures~\ref{fig:fractal} and~\ref{fig:fractalzoom}.

\section{Multiple roots and infinite derivatives}
Of course this method won't help with multiple roots.  We're dividing by derivatives, which go to zero if we approach a multiple root.  It's true that the method \emph{can} work, but if it does it will work more slowly; even so, the sequence $y_k/(y_{k-1}y_{k-2})^2$ will just grow without bound.  I tried this, on $f(x) = (x-2)^2$, with $x_0 = 0.5$, and  this expectation is correct.  The method did converge in ten iterations to the double root $x=2$, with an error of about $10^{-14}$, but note that Newton's method also converges linearly for this problem.

Another place for failure of rootfinding is when any derivative is infinite---that is, at or near derivative singularities.  The presence of up to the \emph{fourth} derivative in the error formula~\eqref{eq:errorformula} suggests that this method will be more sensitive to singularities than Newton's method is.

\section{Where to, next}
The first question to ask is, can this work for nonlinear \emph{systems} of equations?  Even though this method, stripped of its origins, is just a weighted average of two Newton iterations and a secant iteration, all of which can be used for systems, I find that I am dubious.  It's not clear to me what to use for the scalar weights: replace e.g. $y_k^2/(y_k-y_{k-1})^2$ by $\| y_k \|^2/\| y_k - y_{k-1}\|^2$ using some norm?  Which norm?  Perhaps it will work, with the right choice.  But even if it did, for large systems the largest computational effort is usually spent in setting up and solving the linear systems: full Newton iteration is rarely used, when cheaper approximate iterations will get you there faster.  Moreover, the real benefit of the higher speed of convergence of ICI isn't fully felt at double precision or lower---I think it will only matter for quad precision (double-double!) or better.

I believe that the best that can be hoped for from this method is to be used in a computer algebra system for its high-precision scalar solvers, which currently might use other high-order iterations.  It is here, if anywhere, where this method will be useful.


There may be similar two-step iterations of higher order, perhaps starting with ``inverse quintic iteration'' which uses terms up to the second derivative and would have order $(3+\sqrt{21})/2 \approx 3.79$.  I have not yet investigated any such schemes.

In ultra high-precision applications there is sometimes the option to change precision as you go.  This is effective for Newton iteration; only the final function evaluation and derivative evaluation need to be done at the highest precision.  I have not investigated the effects of variable precision on the efficiency of this method.

Then there are more practical matters: what to do in the inevitable case when $y_k = y_{k-1}$? [This happens when the iteration converges, for instance!] Or when $f'(x_k)$ is too small?  It turns out that the details of the extremely practical \texttt{zeroin} method~\cite{Brent:1971:IQI} which combines IQI with bisection and the secant method for increased reliability (this method is implemented in \texttt{fzero} command in Matlab) matter quite a bit.  More improvements on the method can be found in~\cite{wilkins2013modified}, which make the worst-case behaviour better.  Attention needs to be paid to these details.

Then there is another academic matter: in all those tens of thousands of papers on rootfinding, has \emph{nobody} thought of this before?  Really?  Again, I find that I am dubious, although as a two-step method it is technically not one of Schr\"oder's iterations (which keep getting reinvented). And it's not one of Kalantari's methods either~\cite{Kalantari2008polynomial}. But it's such a pretty method---surely someone would have noticed it before, and tried to write about it?  Perhaps it's an exercise in one of the grand old numerical analysis texts, like Rutishauser's (which I don't have a copy of).

There are (literally) an infinite number of iterative methods to choose from; see the references in~\cite{li2019revisiting} for a pointer to the literature, including to~\cite{petkovic2010schroder} and from there to G.~W.~Stewart's translation of Schr\"oder's work from the $19$th century, entitled ``On infinitely many algorithms for solving equations''.  The bibliography~\cite{mcnamee2007numerical} has tens of thousands of entries.

However, searching the web for ``Inverse Cubic Iteration'' (surely a natural name) fails.  We do find papers when searching for names like ``Accelerated Newton's method'', such as~\cite{10.1016/S0096-3003(99)00175-7}; but the ones I have found are each different to ICI (for instance the last-cited paper finds a true third-order method; ICI is not quite third-order but only $1+\sqrt{3}$).  The method of~\cite{kasturiarachi2002leap}, termed ``Leap-Frog Newton'', is quite similar to ICI and consists of an intermediate Newton value followed by a secant step: this costs two function evaluations and a derivative evaluation per step (so is more expensive per step than is ICI which requires only one function evaluation and derivative evaluation per step after the first one) but is genuinely third order.  I found this last paper by searching for ``combining Newton and secant iterations,'' so perhaps I could find papers describing closer matches to ICI, if only I could think of the correct search term. I am currently asking my friends and acquaintances if they know of any.  Do you?

\begin{flushright}
``Six months in the laboratory can save you three days in the library.'' \\
---folklore as told to me by Henning Rasmussen
\end{flushright}
\bibliographystyle{plain}
\bibliography{blendsandsplines}

\begin{figure}[hb]
  \centering
  \includegraphics[width=6cm]{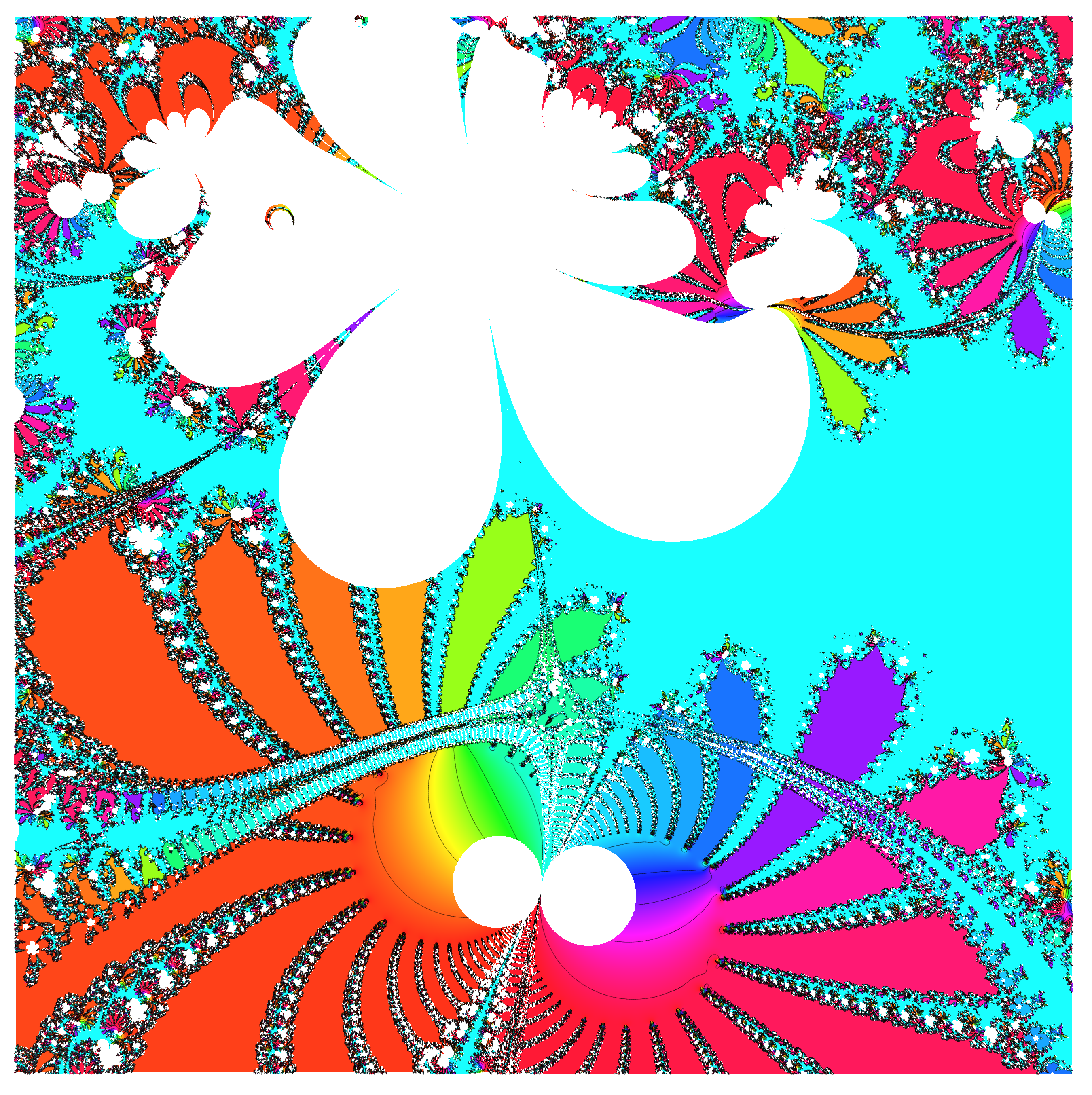}
  \caption{Bonus figure: basins of attraction under ICI for Kepler's equation $z - 0.083\sin(z)-1=0$ on a $1600$ by $1600$ grid $-30.5\le x \le -29.5$, $-17.5\le y \le -16.5$, with tolerance $10^{-8}$ and a maximum of $30$ iterations. The pure white areas represent areas where the iteration encountered a NaN; other colours represent different phases after $30$ iterations or convergence, whichever happened first.  The picture is quite a bit more complicated than the similar picture for Newton's method in Figure 3.21 in~\cite{corless2013graduate}.}\label{fig:Kepler}
\end{figure}

\end{document}